\documentclass{article}      
\usepackage{amssymb}
\usepackage{amsmath}
\usepackage{amscd}
\usepackage{amsthm}  
\usepackage{epsfig}  
\usepackage{subfigure}
\usepackage{wrapfig}
\newtheorem{theorem}{Theorem}

\newtheorem{proposition}[theorem]{Proposition}

\newtheorem{definition}{Definition}
\newcommand{\bz}{\mathbb{Z}}
\newcommand{\br}{\mathbb{R}}

\newcommand{\bh}{\mathbb{H}}

\newcommand{\bs}{\mathbb{S}}
\newcommand{\bl}{\mathbb{L}}
\newcommand{\p}{\partial}
\newcommand{\cc}{\mathcal{C}}
\newcommand{\ce}{\mathcal{E}}

\newcommand{\cb}{\mathcal{B}}

\newcommand{\hk}{\hookrightarrow}

\newcommand{\med}{\medskip}
\newcommand{\la}{\longrightarrow}
\newcommand{\bfl}{\begin{flushleft}}
\newcommand{\efl}{\end{flushleft}}

\newcommand{\eps}{\epsilon}
 
\newcommand{\mtm}{M^{-\tau}}
\newcommand{\tmtm}{\tilde M^{-\tau}}

\newcommand{\ltm}{LM^{-TM}}

\newcommand{\nee}{\nu_\epsilon (e)}
\newcommand{\emk}{Emb (M, \br^k)}
\newcommand{\rek}{\br^k/(\br^k - \nee)}
\newcommand{\sg}{\Sigma}
\newcommand{\rebr}{\br^r / (\br^r - B_1(0))}

\newcommand{\reek}{\br^k/(\br^k - B_\eps (0))}

\newcommand{\smk}{(\bs_M)_k}
\newcommand{\fmb}{F(M, \bs_M)}

\newcommand{\mte}{\mtm (e)}
\newcommand{\vkm}{V_k(\br^{km})}
\newcommand{\nef}{\theta (\phi)}
\newcommand{\mtme}{M^{-\tau}_m(e)}
\newcommand{\mtne}{M^{-\tau}_n(e)}
\newcommand{\mtmne}{M^{-\tau}_{m+n}(e)}

\begin{document}  
%\initfloatingfigs

\title{Multiplicative properties of Atiyah duality}
\author{Ralph L. Cohen   \thanks{The   author was partially supported by a grant from the NSF
}   \\
 Stanford University 
 }
   \date{\today}
 \maketitle

 %\tableofcontents

 \begin{abstract}
Let $M^n$ be a  closed, connected $n$-manifold.  Let $\mtm$ denote the Thom spectrum of its stable normal bundle.    A well known theorem of Atiyah  states that $\mtm$ is homotopy equivalent to the Spanier-Whitehead dual of $M$ with a disjoint basepoint, $M_+$.  This dual can be viewed as the function spectrum, $F(M, S)$, where $S$ is the sphere spectrum.  $F(M, S)$ has the structure of  a commutative, symmetric ring spectrum in the sense of \cite{hss}, \cite{ship} \cite{mayetal}.  In this paper we prove that $\mtm$ also has a natural, geometrically defined,    structure of a commutative, symmetric ring spectrum, in such a way that  the classical duality maps of Alexander, Spanier-Whitehead, and Atiyah define an equivalence of symmetric ring spectra,   $\alpha : \mtm \to F(M, S)$.   
We discuss applications of this to Hochschild cohomology representations of the Chas-Sullivan loop product in the homology of the free loop space of $M$.

\end{abstract}

 \section*{Introduction}  

Throughout this paper $M^n$ will denote a fixed,   closed, connected $n$-manifold.  Let $e : M \hk \br^k$ be an embedding into Euclidean space, and let $\eta_e$ be the normal bundle.  A well known theorem of Atiyah \cite{atiyah} states that the Thom space, $M^{\eta_e}$, is  a $k$- Spanier - Whitehead dual of $M$. One can normalize with respect to $k$ in the following way.  Let $\mtm$ denote the spectrum, $\mtm = \Sigma^{-k}M^{\eta_e}$.  The homotopy type of this spectrum is well defined, in that it is independent of the embedding $e$.  Atiyah's theorem can be restated as saying that there is a homotopy equivalence of spectra,
 
$$
\begin{CD}
\alpha : \mtm @>\simeq >> F(M, S)
\end{CD}
$$
where $S$ is the sphere spectrum, and  $F(M, S)$ is the spectrum whose $k^{th}$ space is the unbased mapping space, $F(M, S^k)$.  

Recently, symmetric monoidal categories of spectra have been developed
(\cite{maygroup}, \cite{hss}, \cite{mayetal}).  In these categories, the dual  $F(M, S)$ has the structure  of a commutative ring spectrum.  The goal of this paper is to prove that the Thom spectrum $\mtm$ also has such a multiplicative structure via a natural geometric construction.  Moreover we will show that the classical duality maps of Alexander, Spanier-Whitehead, and Atiyah define an equivalence of ring spectra,
$\alpha : \mtm \xrightarrow{\simeq} F(M, S)$.   

For our purposes the most convenient category of spectra will be the category of symmetric spectra of \cite{hss}, \cite{hovey} \cite{ship}.     We will actually  work in the topological  (as opposed to simplicial) category of symmetric spectra as developed in \cite{mayetal}.  The following  is the main result of  this paper.

\begin{theorem}\label{one}
The Thom spectrum $\mtm$ has the structure of a a commutative, symmetric ring  spectrum with  unit.  Furthermore, the duality map
$$
\alpha : \mtm \to F(M, S)
$$
is a $\pi_*$-equivalence of symmetric ring spectra.  That is,  $\alpha$ is a ring map that induces an isomorphism on stable homotopy groups. 
 \end{theorem}

  The study of the multiplicative properties of $\mtm$ and the Atiyah duality map were motivated by its use in string topology \cite{cohenjones}, and the question of whether the Chas-Sullivan string topology operations defined on the homology of the loop space, $H_*(LM)$,  are homotopy invariants \cite{chsull}.  An important step in this investigation is the relationship between
   the Chas-Sullivan  loop homology  algebra $\bh_*(LM)$ and the Hochschild cohomology $H^*(C^*(M), C^*(M))$.  We will use the multiplicative properties of Atiyah duality proven here to fill in details of the argument given in \cite{cohenjones} showing that these algebras are isomorphic, when $M$ is simply connected.  
   
   While the algebra $H^*(C^*(M), C^*(M))$ is clearly a homotopy invariant of $M$, the ring isomorphism between $\bh_*(LM)$ and $H^*(C^*(M), C^*(M))$ uses embeddings of manifolds and the Thom-Pontrjagin construction, and hence depends a priori on the smooth structure of $M$.  The question of the homotopy invariance of the string topology operations will be the subject of a future paper by the author, J. Klein, and D. Sullivan \cite{cks}.  
  
  The organization of this paper is as follows.  In section 1 we describe the geometric construction of the  symmetric spectrum $\mtm$ using Thom spaces parameterized by spaces of embeddings and tubular neighborhoods.  We then describe the ring structure.  It turns out that this spectrum does not naturally have a unit, essentially because there is no canonical choice of embedding.  We address this issue in section 2  by showing that by choosing
  a fixed embedding  $e: M \hk \br^k$ one can construct an equivalent symmetric ring spectrum $\mtm (e)$ that does admit a unit.  This ring spectrum
  will be in the sense of \cite{hovey} where the underlying sphere spectrum is made out of iterated smash products of $S^k$, rather than $S^1$ as in \cite{hss}.   We then complete the proof of theorem \ref{one} in section 2.     We end by giving  the applications to Hochschild cohomology mentioned above. 
  
  The author would like to thank G. Carlsson, B. Dundas, J.P. May, D. Sullivan, C. Schlichtkrull, and T. Tradler for many helpful
  conversations and correspondence related to this work.

\section{The Thom spectrum $\mtm$}

Let $M^n$ be a closed $n$-manifold.  The homotopy type of the spectrum $\mtm$ is the desuspension of the Thom space $\Sigma^{-k}M^{\eta_e}$, where $\eta_e$ is the normal bundle of an embedding, $e : M^n \hk \br^k$.  In order to describe $\mtm$ as a symmetric ring spectrum we need to deal with the ambiguities in this description, such as choices embeddings,  tubular neighborhoods, etc.  We do this by keeping track of all such choices in our definition of $\mtm$. 

Let $e: M \hk \br^k$ be an embedding.  Recall that the tubular neighborhood theorem says that for $\eps$ small enough, there is an open set $\nee$ containing $e(M)$ and a retraction $r : \nee \to e(M)$ satisfying the following properties:
\begin{enumerate}\label{tube}
\item $\|x-r(x)\| \leq \|x-e(y)\| $\, \, for any  $y \in M$,  \, \, with equality holding if and only if $r(x) = e(y)$.
\item For all $y \in M$, $r^{-1}(e(y))$ is a ball of radius $\eps$ in the affine space $e(y) + T_{e(y)}e(M)^\perp$, with center at $e(y)$. 
\item The closure $\bar \nu_\eps (e)$ is a smooth manifold with boundary.
\end{enumerate}
We remark that the retraction $r$ is completely determined by properties (1) and (2).  See \cite{madtorn} for this version of the tubular neighborhood theorem.
Let $L_e >0 $ be  the minimum  of $1$ and the least upper bound of those $\eps$ satisfying this theorem.

For $k >0$, define the space $\tmtm_k$ as follows.
 \begin{equation}\label{tmtm}
 \tmtm_k = \begin{cases} point \quad \text{if} \quad \emk = \emptyset ,\\
 \{(e, \eps, x) \, : e: M \hk \br^k  \, \text{ is an embedding,}  \,  0<\eps< L_e, \\ \text{and} \,   x \in \rek \}  \quad  \text{if} \quad  \emk \neq \emptyset
 \end{cases}
 \end{equation}
 
Let $\ce_k = \{(e, \eps) : \, e \in Emb(M, \br^k)$, and $\eps \in (0, L_e)\}$.  This space is topologized as a subspace of $Emb(M, \br^k) \times (0, \infty)$.   Notice that  $\tmtm_k$ can be topologized as the total space of a fiber bundle  
 \begin{equation}\label{fiber}
 p : \tmtm_k \to \ce_k
 \end{equation}
  whose fiber over $(e, \eps)$ is  $\rek$, which by the tubular neighborhood theorem is  homeomorphic to    the Thom space of the normal bundle, $\eta_e$. Notice also that there is a canonical section of this bundle
 $$
 \sigma_\infty : \ce_k \to \tmtm_k
 $$
 given by $(e, \eps) \to (e, \eps, \infty)$, where $\infty$ is the basepoint in $\rek$.   
 We now define
 \begin{equation}\label{mtm}
 \mtm_k = \tmtm_k/(\sigma_\infty (\ce_k).
 \end{equation}
 
 \med
 
   A point in $\mtm_k$ is an embedding with tubular neighborhood, together with a point in the Thom space $\rek$.  
 
 \med
 Now let $(e, \eps) \in \ce_k$ be fixed.   Let
 $
 \phi : M^{\eta_e} \to \rek
 $
 be the homeomorphism given by the tubular neighborhood theorem.  This defines an inclusion
 $  j_e: M^{\eta_e} \xrightarrow{\phi}   \rek \xrightarrow{x \to (e, \eps, x)}   \mtm_k.
 $ 
 
 The following states that through a range of dimensions that increases with $k$, $\mtm_k$ is just the Thom space $M^{\eta_e}$.
 
 \begin{proposition}\label{range}
 The inclusion $j_e : M^{\eta_e} \to \mtm_k$ induces an isomorphism in homotopy groups through dimension $\frac{k}{2}-n-2$.
 \end{proposition}
 
 \begin{proof}  This follows from fibration (\ref{fiber})  and the Whitney embedding theorem which states that $\emk$ is $(\frac{k}{2}-n-1)$-connected,  and therefore so is  $\ce_k$.  
 \end{proof}
 
 \med
 Notice that the permutation action of the symmetric group $\Sigma_k$  on $\br^k$ induces a $\Sigma_k$-action on $\ce_k$, and on $\mtm_k$.   This gives $\{\mtm_k\}$ the structure of a symmetric sequence in the sense of \cite{hss}.
 
 \med
 Let $S^k$  be the unit sphere,        $S^k = \br^k / (\br^k - B_1(0))$, where
 $B_1(0)$ is the unit ball around the origin.  The $\sg_k$-action on $\br^k$ induces a $\sg_k$-action on $S^k$.    
 
 We  define structure maps  
  \begin{align}\label{structure}
 s_{r,k} : S^r \wedge \mtm_k  &\la \mtm_{r+k} \\
  t \wedge (e, \eps, x) &\to (0 \times e, \eps, t \wedge x) \notag \\
\end{align}
 where $t\wedge x \in \br^{r+k}/( \br^{r+k}- \nu_{\eps}(0\times e))$ is the image of the smash product of $t$ and $x$ under the projection map
 $$
 \rebr \wedge \rek = \br^{r+k}/( \br^{r+k}- (B_1(0)\times \nu_{\eps}( e))) \la \br^{r+k}/( \br^{r+k}- \nu_{\eps}(0\times e)).
 $$
 Clearly the maps $s_{r,k}$ are $\sg_r \times \sg_k$-equivariant.  They also satisfy the appropriate compatibility conditions to allow us  to    make the following definition.
 
 \med
 \begin{definition}\label{mtmdef}
 The symmetric spectrum $\mtm$ is defined by the collection $\{\mtm_k, s_{r,k} : S^r \wedge \mtm_k \to \mtm_{r+k}\}$.
 \end{definition}
 
 \med
 Notice that by proposition \ref{range},  the spectrum $\mtm$ has the right homotopy type.  That is,  given an embedding $e : M \hk \br^k$   then
 there is an equivalence of spectra, $ \sg^\infty M^{\eta_e} \to \sg^k \mtm$. 
 
  \med
We now   examine the multiplicative structure of $\mtm$.  Recall from \cite{ship}  \cite{mayetal},  that a (commutative) symmetric ring spectrum is a (commutative) monoid in the symmetric monoidal category of symmetric spectra.  The spectrum $\mtm$ as we have defined it will have an associative, commutative product structure, but it does not have a unit.  This is essentially because there is no canonical choice of embedding  of $M$ in $\br^k$. We will deal with issues regarding the unit in the next section.  

\med
\begin{theorem}\label{ring}$\mtm$ is a commutative symmetric ring spectrum without unit.  That is, it is a commutative monoid  without unit,(i.e ``commutative semigroup")  in the category of symmetric spectra.
\end{theorem}

\begin{proof}  We define associative pairings 
$$
\mu_{k,r} : \mtm_k \wedge \mtm_r \la \mtm_{k+r}
$$
by the formula
 \begin{equation}\label{mukr}
\mu_{k,r} ((e_1, \eps_1, x_1), (e_2, \eps_2, x_2)) = (e_1 \times e_2, \, \eps_{1,2}, \, x_1\wedge x_2)
\end{equation}
where  $e_1 \times e_2$ is the embedding 
$  M \xrightarrow{\Delta}  M\times M \xrightarrow{e_1 \times e_2}  \br^k \times \br^r  $, 
$\eps_{1,2} = min\{\eps_1, \eps_2, L_{e_1\times e_2}\}$, and $x_1\wedge x_2$ is the image of the smash product under the obvious projection map,
 \begin{multline}
\br^k/(\br^k - \nu_{\eps_1}(e_1)) \, \wedge\,  \br^r / (\br^r-\nu_{\eps_2}(e_2)) = \br^{k+r}/(\br^{k+r}-(\nu_{\eps_1}(e_1) \times\nu_{\eps_2}(e_2)))  \\
 \la \br^{k+r}/(\br^{k+r}-\nu_{\eps_{1.2}}(e_1\times e_2)).
 \end{multline}
  
  The map  $\mu_{k,r}$ is clearly $\sg_k \times \sg_r$-equivariant, and the collection $\{\mu_{k,r}\}$  is associative.     Furthermore a check of definition of the structure maps (\ref{structure}) shows that these pairings respect the structure maps, and so define a pairing on the external tensor product of spectra (\cite{hss} \cite{mayetal})
 $$
 \mu : \mtm \otimes \mtm \la \mtm.
 $$
   In order to show that the pairing $\mu$ defines a ring structure, we need to show that it factors through the smash product,
 $\mtm \wedge \mtm$, defined to be the coequalizer of 
 $$
 r \otimes 1 \quad \text{and} \quad 1 \otimes \ell \quad : \mtm \otimes S \otimes \mtm \la \mtm \otimes \mtm
 $$
 where $\ell : S \otimes \mtm \to \mtm$ is the map induced by the structure map $s_{r,k} : S^r \wedge \mtm_k \to \mtm_{r+k}$, and $r : \mtm \otimes S \to \mtm$ is induced by the maps
 $
  r_{p,q} : \mtm_p \wedge S^q  \xrightarrow{switch}   S^q \wedge \mtm_p \xrightarrow{s_{q,p}}  \mtm_{p+q} \xrightarrow{\tau_{p,q}} \mtm_{p+q},$ where $\tau_{p,q} \in \sg_{p+q}$ is the  shuffle permutation,
 $$
 \tau_{p,q} (i) = \begin{cases}  q+i, \quad i \leq i \leq p \\
 i-p,  \quad  p+1 \leq i \leq p+q.  \end{cases}
 $$
This fact can now be immediately checked   from the definitions. 
 To show that the ring structure is commutative, we need to verify that the following diagrams commute, for every $p$ and $q$.
 
 $$
 \begin{CD}
 \mtm_p \wedge \mtm_q  @>\mu_{p,q} >>  \mtm_{p+q} \\
 @Vswitch VV     @VV\tau_{p,q} V \\
 \mtm_q \wedge \mtm_p   @>>\mu_{q,p} >  \mtm_{p+q}
 \end{CD}
 $$
    Again, this can be seen by a quick check of the definitions.
\end{proof}

\med

We now use this ring spectrum to study  the Atiyah duality map.    
We begin by recalling the classical map inducing 
ity in homology. 
Namely, let $e : M \hk \br^k$ be an embedding with tubular neighborhood $\nee$.  Consider the map
\begin{align}\label{alex}
\alpha : \left( \br^k - \nee\right) \times M & \to \br^k - B_\eps (0) \simeq S^{k-1} \notag \\
(v, y)  &\la v - e(y).
\end{align}
This map induces the Alexander duality isomorphism
$$
\begin{CD}
 \tilde H_q(\br^k - e(M)) \cong \tilde H_q(\br^k - \nee) @>\cong >> \tilde  H^{k-q-1}(M).
\end{CD}
$$

Atiyah duality \cite{atiyah}  is induced by   the same map:
\begin{align}
M^{\eta_e} \wedge M_+ \cong (\br^k \times M )/\left( (\br^k - \nee) \times M \right)  &\la \br^k /( \br^k - B_\eps (0)) \quad \simeq S^k  \notag \\
(v, y)   &\la v-e(y).
\end{align}
The adjoint of this map gives a map
$ 
\alpha : M^{\eta_e}  \la F(M, S^k)
$
which defines, up to homotopy, the Atiyah duality homotopy equivalence,
\begin{equation}\label{atdual}
\alpha : \mtm \la F(M, S).
\end{equation}

To see this map as a map of symmetric ring spectra, we introduce a function spectrum
$\fmb$,  where $\bs_M$ is a symmetric ring spectrum (without unit), equivalent to the sphere spectrum $S$.  

\med
   Define
\begin{equation}\label{tildesm}
(\tilde \bs_M)_k = 
 \{(e, \eps,   t) \; : \; (e,\eps) \in \ce_k,  \; \text{and} \; t \in \reek\}.
 \end{equation}
 $(\tilde \bs_M)_k$ is topologized as a fiber bundle over $\ce  $, with fiber over $(e, \eps, v)$ equal to $\reek$.  

The bundle $p:(\tilde \bs_M)_k\to \ce_k  $ has a canonical section $\sigma_\infty  (e, \eps) = (e, \eps,  \infty)$.  We   then define the space
$\smk$ to be the cofiber,
\begin{equation}\label{smk}
\smk = (\tilde \bs_M)_k/\sigma_\infty (\ce_k ).
\end{equation}

\med

$\smk$ has the $\sg_k$-action induced by the permutation action on $\br^k$.  We  define structure maps
$$
\beta_{r,k} : S^r \wedge \smk \to (\bs_M)_{r+k}
$$
by $s \wedge (e, \eps;  t) \la (0 \times e, \eps,  s\wedge t)$, where $0 \times e$ is the embedding $y \to (0, e(y)) \in \br^r \times \br^k$.    This gives $\bs_M$ the structure of a symmetric spectrum.
The ring structure on $\bs_M$ is given by the pairings
$$
m_{r,s} : (\bs_M)_{r} \wedge (\bs_M)_s \to (\bs_M)_{r+s}
$$
defined by 
\begin{equation}\label{smring}
(e_1, \eps_1,  t_1) \wedge (e_2, \eps_2,  t_2) \la (e_1\times e_2, \eps_{1,2},   t_1\wedge t_2).
\end{equation}
 It is straightforward to check that these pairings give $\bs_M$ the structure of a commutative, symmetric ring spectrum (without unit).
Furthermore, there is an obvious map $\bs_M \to S$ which is  a $\pi_*$-equivalence of symmetric ring spectra.  

We now consider the symmetric ring spectrum $\fmb$.  (Its ring structure is induced by that of $\bs_M$.)  Notice that the Atiyah duality map  defines a map
\begin{align}\label{define}
\mtm_k \wedge M_+   &\la  (\bs_M)_k \\
(e, \eps, x) \wedge y  &\to (e, \eps, x - e(y))
\end{align}
where, as above, $x \in \rek$, and so $x-e(y) \in \reek$.  A quick check of definitions verifies that the adjoint of this map defines a map of symmetric ring spectra
$$
\alpha : \mtm \la \fmb
$$
realizing the Atiyah homotopy equivalence.  

\med
Unfortunately there are no good unital properties of this equivalence of ring spectra.
This problem will be dealt with in the next section.

\section{Unital symmetric ring spectra and a proof of theorem \ref{one}}

\med
In this section we use a fixed embedding $e : M \hk \br^k$ to  modify the definitions of $\mtm$ and $F(M, S)$, in order to obtain symmetric ring spectra with units. We then use these spectra   to prove theorem \ref{one}.
To do this we use the machinery of \cite{hovey}, which defines symmetric spectra in general   model categories. 
If $\cc$ is a symmetric, monoidal, model category, and $K$ is a cofibrant object of $\cc$, then Hovey defined a category $Sp^\Sigma (\cc, K)$ of symmetric spectra.  These are symmetric sequences $\{X_n\}$ together with $ \Sigma_n \times \Sigma_m$-equivariant structure maps, $\eps_{n,m}: X_n \otimes K^{\otimes m} \to X_{n+m}$.  In the standard setting considered above, $\cc$ is the category of based topological spaces, and $K = S^1$.  For our purposes below, we will continue to let $\cc$ be the category of based topological spaces, but we now let $K = S^k$, where $k$ is the ambient dimension of our fixed embedding, $e : M \hk \br^k$, and $S^k = \br^k \cup \infty$ is the one point compactification. 

We now define a symmetric ring spectrum $\mte$  in this category.    In the definition of the spaces $\mtm_n$, we considered all possible embeddings
of $M$ in $\br^n$ together with tubular neighborhoods.  In our present situation, the fact that we have a fixed embedding $e : M \to \br^k$, allows us to restrict the general embeddings we consider to those of the form
$\phi \circ e$, where $\phi : \br^k \to \br^{kn}$ is a linear, isometric embedding.  

Let $0 < \eps < \frac{1}{2}L_e$ be fixed, and $\nee \subset \br^k$ be the corresponding tubular neighborhood.  Let $\vkm$
 be the Stiefel manifold of $k$-frames in $\br^{km}$, or equivalently, linear isometric embeddings of $\br^k$ in $\br^{km}$.  This Stiefel manifold has a natural action of the symmetric group $\Sigma_m$ given by permuting
 the factors of $(\br^k)^m = \br^{km}$.    $\vkm$ has a   basepoint given by
 the $m$-fold diagonal embedding, $\Delta_m : \br^k \to (\br^k)^m$.

 For $\phi \in \vkm$, define the open neighborhood of the embedding $\phi \circ e : M \hk \br^{km}$
 $$
 \theta (\phi) = \{z \in \br^{km} \ : \ \exists x \in \nee \, \text{with} \, |z- \phi (x)| < \eps \}.
 $$ Notice that $\theta (\phi ) \subset \nu_{2\eps}(\phi \circ e)$ and is a regular neighborhood of the embedding $\phi \circ e (M)$. 
    Define the spaces
 $$
 \tilde M^{-\tau}_m (e) =\{(\phi, x): \phi \in \vkm, \, \text{and} \, x \in \br^{km}/(\br^{km} - \nef) \}.
 $$
 This is topologized as a fiber bundle over $\vkm$, with fiber the  appropriate Thom space.   $  \tilde M^{-\tau}_m (e)$
 has an  action of $\Sigma_m$ induced by the permutation action on $\br^{km}$, making this bundle $\Sigma_m$-equivariant.  As before, let $\infty \in  \br^k/(\br^k - \nef) \}$ denote the basepoint.  We can then define the quotient
 
 \begin{equation}\label{mte}
 \mtme =   \tilde M^{-\tau}_m (e) /\{(\phi, \infty) : \phi \in \vkm\}.
 \end{equation}
 
 \med
 This now defines is a symmetric sequence in the  category $\cc$ of based topological spaces. (For this we take $M^\tau_0(e) = S^0$.)   
 
 This symmetric sequence is multiplicative, in the sense that we have associative pairings
\begin{align} 
 \mu_{m,n} : \mtme \wedge \mtne &\to \mtmne \notag \\
 (\phi_1, x_1) \wedge (\phi_2, x_2) &\to (\phi_1 \times \phi_2, x_1 \wedge x_2) \notag
 \end{align}
  where  $\phi_1 \times \phi_2 \in V_k(\br^{km}\times \br^{kn})$ is the product of the linear embeddings $\phi_1$ and $\phi_2$.  The element $x_1 \wedge x_2$ is the image of the smash product under the projection map
  $$
  \br^{km} \times \br^{kn}/( \br^{km} \times \br^{kn} - \theta (\phi_1) \times \theta (\phi_2)) \to
  \br^{k(m+n)}/ (  \br^{k(m+n)} - \theta(\phi_1\times \phi_2) ).
  $$
  Notice that these pairings are also commutative in the sense that   $$\mu_{n,m} \circ switch  = \tau_{n,m}\circ \mu_{m,n}: \mtme \wedge \mtne  \to \mtmne$$ where $\tau_{m,n}$ is the  shuffle permutation
  described in the last section. 
  
 Now consider the symmetric sequence $S(k)$, defined by $S(k)_m = (S^k)^{\otimes m}$.  Since tensor product
 in the category $\cc$ is given by smash product, these spaces are just spheres.   This sequence has, by its definition, an associative product structure, which is commutative in the above sense.

\med
We now construct a multiplicative map of symmetric sequences,
\begin{equation}\label{newunit}
u : S(k) \to \mte
\end{equation}
which will serve as the unit.  Define
\begin{align}
u : S^{km} &\to \mtme \notag \\
t &\la (\Delta_m, c (t))\notag
\end{align}
where Thom collapse map $c : S^{km} = \br^{km} \cup \infty  \to \br^{km}/(\br^{km}-\theta (\Delta_m))$ is the projection.  An easy check of the definitions verifies that 
  $u : S(k) \to \mte$ is a multiplicative map of symmetric sequences.
 This in particular defines a bimodule structure of 
 $\mte$  over $ S(k)$.     These definitions yield the following.
 
 \med
 \begin{proposition}
 The bimodule structure of the symmetric sequence $\mte$ over  $S(k)$,   and the pairings $\mu_{m,n}: \mtme \wedge \mtne \to \mtmne$ define a symmetric ring spectrum in the category
 $Sp^\Sigma (\cc, S^k)$.  $u : S(k) \to \mte$ is the unit in this ring structure.
 \end{proposition}
 
 \med
 
 Notice that since the connectivity of the Stiefel manifold $\vkm$ goes up with $m$, the homotopy type of the spectrum $\mte$ is the same as the homotopy type of the desuspension of the Thom spectrum of the normal bundle, $\Sigma^{-k}M^{\eta_e}$.  This in turn is the same as the homotopy type of the spectrum $\mtm$ defined in the last section.

 \bf Remark.  \rm  The definition of $\mte$ depended on a fixed choice of $\eps < \frac{1}{2}L_e$.  Normally we will suppress the choice of $\eps$ in the notation.  However when we want to make note of it, we will denote the resulting spectrum by $M^{\tau}(e, \eps)$.  Notice that for $\eps^\prime < \eps$, there is a natural projection
 map $M^{\tau}(e, \eps) \to M^{\tau}(e, \eps^\prime)$ which is an equivalence of symmetric ring spectra.
  
  \med
 We now define a variant of the Spanier-Whitehead dual of $M$ in the category $Sp^\Sigma (\cc, S^k)$.
  For $\phi \in \vkm$, define a fiber bundle over the  tubular neighborhood,
 $$
\pi_\phi :  \cb_\phi \la   \nu_\eps (e)
 $$
 where the fiber over $x \in   \nu_\eps (e)$, is the sphere, $\bar B_{\eps}(x)/\p\bar B_{\eps}(x)$.  Here $ B_{\eps}(x)$ is the ball in $\br^{km}$ of radius $\eps$ around $\phi(x)$, and  $\bar B_{\eps}(x)$ is its closure. The permutation action on $\br^{km}$ induces a $\Sigma_m$ action on this bundle,   where $\Sigma_m$ acts trivially on the base.     
 
 Let $\rho : ([0,\eps), [0, \frac{\eps}{2}) ) \to( [0, +\infty), [0, \frac{\eps}{2}))$ be an explict diffeomorphism of the half open
 interval $[0, \eps) $ with $[0, +\infty)$ which fixes $[0,\frac{\eps}{2})$. This map defines radial expansion diffeomorphisms,
 $\rho_x : (B_\eps (x), B_{\frac{\eps}{2}}(x)) \to ( \br^{km},  B_{\frac{\eps}{2}}(x))$, for each $x \in \nee$.  This defines a trivialization of the sphere bundle,
 
    \begin{equation} \label{trivial}
 \rho : \cb_\phi \to \nee \times  S^{mk}  
   \end{equation}
 
 Let $\Gamma(\cb_\phi)$ be the space of sections of $\cb_\phi $.  This trivialization induces
 a $\Sigma_m$-equivariant homeomorphism $\rho : \Gamma(\cb_\phi) \xrightarrow{\cong} F(\nee, S^{km})$. Restriction defines an equivariant homotopy equivalence, $F(\nee, S^{km}) \xrightarrow{\simeq} F(M, S^{km}).$
 
 \med
 Now define 
 \begin{equation}\label{fm}
 F_m(e) = \{(\phi, \sigma) : \, \phi \in \vkm, \, \sigma \in \Gamma(\cb_\phi) \}/\{(\phi, \sigma_\infty)\},
 \end{equation}
 where we are dividing out by all pairs $(\phi, \sigma_\infty)$, where $\sigma_\infty$ is the constant section at the basepoint.  Notice that the trivialization $\rho$ defines an $\Sigma_m$-equivariant homeomorphism,
\begin{equation}\label{fmtriv}
 \rho : F_m (e) \xrightarrow{\cong} \vkm_+ \wedge F(\nee, S^{mk}),
 \end{equation}
 which in turn is equivariantly homotopy equivalent to $ \vkm_+ \wedge F(M, S^{mk})$.
  Let $ \bar \rho: F_m (e) \to F(M, S^{mk})$ be the composite $F_m (e) \xrightarrow{\rho} \vkm_+ \wedge F(\nee, S^{mk}) \cong \vkm_+ \wedge F(M, S^{mk}) \xrightarrow{project} 
 F(M, S^{mk})$. 
 
 \med
 We now have constructed a symmetric sequence $F(e) = \{F_m(e)\}$.  This sequence is multiplicative, in that it has an associative multiplication,
\begin{align}
 \mu_{m,n} : F_m(e) \wedge F_n(e) &\to F_{m+n}(e) \\
 (\phi_1, \sigma_1) \wedge (\phi_2, \sigma_2) &\to (\phi_1\times \phi_2,  \sigma_{1,2}) \notag
 \end{align}
 where $\sigma_{1,2} (x) = \sigma_1(x)\wedge \sigma_2(x)$.  This multiplication is also commutative in the sense that   $\mu_{n,m} \circ switch  = \tau_{n,m}\circ \mu_{m,n}: F_m(e) \wedge F_n(e) \to F_{m+n}(e).$ 
 
 Define the Atiyah duality map  $\alpha : \mtm \to F(e)$,   as follows.   
\begin{align}\label{alpha}
 \alpha : \mtme &\la F_m(e) \notag\\
 (\phi, x) &\to (\phi, \sigma_x)
 \end{align}
 where $\sigma_x (y)$  is defined to be the basepoint if $x = \infty$, or if $|x - \phi (y)|> \eps$.  Otherwise $\sigma_x (y)  = [x] \in \bar B_{\eps}(y)/\p \bar B_{\eps}(y)$.   
  One immediately checks that this is a map of symmetric sequences (i.e is equivariant at each stage), and preserves the multiplicative structure.  The composition with the unit,
 $$
 S(k) \xrightarrow{u} \mtm \xrightarrow{\alpha} F(e)
 $$
defines a bialgebra structure of $F(e)$ over $S(k)$, giving
 $F(e)$ the structure of a symmetric ring spectrum in $Sp^\Sigma (\cc, S^k)$.  This composition defines a unit for $F(e)$. Thus
  the Atiyah duality map
 $
 \alpha : \mtme \to F(e)$
is a map of (unital) symmetric ring spectra.    

Furthermore, with this ring spectrum structure, the maps    $\bar \rho : F_m(e) \to F(M, S^{km})$ induce an equivalence of symmetric ring spectra, 
$$
\bar \rho: F(e) \xrightarrow{\simeq} F(M,S).
$$
By comparing formulas (\ref{trivial}) and (\ref{alpha}), the composition of $\bar \rho$ with $\alpha$ is precisely the classical Atiyah duality equivalence
described in the last section.  Since these are both equivalences of symmetric ring spectra, so is their composition, which  proves theorem \ref{one}.

\med
We now apply these results to Hochschild homology questions.  Since $\alpha : \mtme \to F
(e)$ is an equivalence of symmetric ring spectra, then they induce equivalences of their
topological Hochschild homologies  and cohomologies:
$$
THH_*(\mtme) \simeq THH_*(F(e))\simeq THH_*(F(M, S)) 
$$
$$
THH^*(\mtme) \simeq THH^*(F(e))\simeq THH^*(F(M,S)).  
$$
Furthermore, by results of McClure and Smith \cite{mcsm}, the topological Hochschild cohomology
of a symmetric ring spectrum has the structure of an algebra over the little disk operad $\cc_2$,
and hence the equivalence above  preserves this multiplicative structure.    This has the following consequence when we apply 
singular chains functor, $C_*(-)$.

Since the evaluation map induces a  chain homotopy equivalence of differential graded algebras,
$
ev_* : C_*(F(M,S)) \xrightarrow{\cong} C^*(M),
$
we have   an isomorphism of Hochschild cohomologies as  algebras over the  operad, $H_* (\cc_2)$,
\begin{equation}\label{isoone}
H^*(C_*(\mtm), C_*(\mtm)) \cong H^*(C^*(M), C^*(M)).
\end{equation}
(We are surpressing the embedding $e$ from this notation, since the chain homotopy type of $C_*(\mtme)$ is clearly independent of the embedding
$e$.  )

\med
The differential graded algebra $C_*(\mte)$ is also a bimodule over $C^*(\nee)$.  This is seen as follows.  Consider the diagonal
maps
$$
\Delta_r : \mte \to \mte \wedge \nu_{2\eps}(e)_+ \quad \text{and} \quad \Delta_{\ell} : \mte \to \nu_{2\eps}(e)_+\wedge \mte
$$
defined by sending $(\phi, x)$ to $x\wedge x_1$ and $x_1 \wedge x$ respectively.  Here,   if $x \in \theta(\phi)   \subset \br^{km}$, then $x = \phi (x_1) \oplus x_2$ is the unique decomposition where $x_1 \in \nu_{2\eps}(e) $ and $x_2 \in \phi (\br^k)^\perp$. 
By applying chains and composing with the evaluation map, $C^*(\nu_{2\eps}(e)) \otimes C_*(\nu_{2\eps}(e)) \to \bz$,
we have a bimodule structure of $C_*(\mte)$ over $C^*(\nu_{2\eps}(e))$.   Now let $q: \nu_{2\eps}(e) \to \nee$ be a fixed retraction.  Then composing with $q^* : C^*(\nee) \to C^*(\nu_{2\eps}(e))$ defines the bimodule structure over $C^*(\nee)$.

This bimodule structure can be realized on the spectrum level as follows.  Define the map
\begin{align}\label{module}
r_{m,n} : \mtme \wedge F_n(e)  &\la M^{-\tau}_{m+n}(e) \notag\\
(\phi_1, x) \wedge (\phi_2, \sigma) &\la (\phi_1 \times \phi_2, x \wedge u(\bar \rho (\sigma(q(x_1)))) 
\end{align}
where $u$ is the unit, and $\bar \rho$ is the trivialization described above.      There is also a similar map
$$
\ell_{n,m} : F_n(e) \wedge \mtme \to M^{-\tau}_{m+n}(e),
$$ which together define a bimodule structure of $C_*(\mtme)$ over $C_*(F(e))$.  With respect to the   chain equivalence   $
 C_*(F(e)) \xrightarrow{\rho}  C_*(F(\nee, S))\xrightarrow{ev}  C^*(\nee),
 $
 this bimodule structure extends  the bimodule structure of $C_*(\mte)$
 over $C^*(\nee)$   above. This structure respects the Atiyah duality map $\alpha$ in the following sense.
 Let $p: M^{-\tau}(e,\eps) \xrightarrow{\simeq} M^{-\tau}(e,\frac{\eps}{2})$ be the equivalence of ring spectra
 defined by projecting the spectrum defined with a fixed $\eps$ to the one defined with $\frac{\eps}{2}$.  A check of definitions verifies that  the compositions $C_*(M^{-\tau}(e, \eps)) \otimes   C_*(M^{-\tau}(e, \eps)) \to  C_*(M^{-\tau}(e,\frac{\eps}{2})) $ defined by $p \circ r \circ (1 \wedge \alpha)$ and $p \circ \ell \circ (\alpha \wedge 1)$
 are both equal to $p \circ \mu : C_*( M^{-\tau}(e,\eps)) \otimes  C_*(M^{-\tau}(e,\eps)) \to C_*(M^{-\tau}(e,\frac{\eps}{2}))$.

By considering this module action applied to the unit in $C_*(\mte)$, one has a map
 $$
 \beta_* : C_*(F(e)) \to C_*(\mte),
 $$
 which by a check of definitions (\ref{newunit}), (\ref{alpha}), and (\ref{module}) is immediately
 seen to be a ring homomorphism, and by the above observation is    a left inverse to the Atiyah duality map,  $ \alpha_* : C_*(\mte) \to C_*(F(e)) $.    
 As a consequence of these observations, we can draw the following conclusion.
 
 \begin{proposition}\label{hoch} The Hochschild cohomology $H^*(C^*(M),C_*(\mte)) \cong H^*(C^*(\nee), C_*(\mte)) \cong H^*(C_*(F(e)), C_*(\mte))$ is an algebra isomorphic (as algebras) via the Atiyah duality map to $H^*(C_*(\mte), C_*(\mte))$, which  in turn is isomorphic as algebras  to
 $H^*(C^*(M), C^*(M))$.
 \end{proposition}
 
 We use this result to fill in details of the proof of the following theorem \cite{cohenjones}.
 
 \begin{theorem}Let $M$ be a simply connected closed manifold, $LM$ its free loop space,
 and $\bh_*(LM)= H_{*+d}(LM)$ the loop homology algebra of Chas-Sullivan \cite{chsull}.  Then there is an isomorphism of algebras,
  $$
 H^*(C^*(M), C^*(M)) \cong \bh_*(LM).
 $$
 \end{theorem}

In \cite{cohenjones}, the authors described a ring spectrum $\ltm$ defined as the Thom spectrum of the pull-back of the virtual bundle $-TM \to M$ over the loop space $LM$ via the map
$
e: LM \to M
$
which evaluates a loop at the basepoint $1 \in S^1$.  It was shown that upon application of the Thom isomorphism in homology, $\tau_* : H_*(\ltm) \to H_{*-d}(LM),$  the ring structure of $\ltm$ corresponds to the Chas-Sullivan product.  Furthermore, for $M$ simply connected,  a cosimplicial spectrum model $(\bl_M)_*$  for $\ltm$ was described in section 3 of \cite{cohenjones} whose spectrum of  $k$-cosimplices of $\bl_M$ was given by 
$
(\bl_M)_k =  
  = \mtm  \wedge (M^k)_+.
$  For the purposes of this note we replace this by the homotopy equivalent spectrum of cosimplices, $(\bl_M)_k = \mte \wedge (\nee_+)^{(k)}.$  
   As in \cite{cohenjones},  the   coface and codegeneracy maps are defined in terms of the diagonal maps described above, composed with the retraction $q : \nu_{2\eps}(e) \to \nee$.  A specific homotopy equivalence  was given in \cite{cohenjones},
   $$
   f: \ltm \xrightarrow{\simeq} Tot(\bl_M).
   $$
   The chains  of the $k$-cosimplices  are   given by
   $$
C_*(\mtm) \otimes C_*(M)^{\otimes k} \cong Hom (C^*(M)^{\otimes k}, C_*(\mtm)) \cong  Hom (C^*(\nee)^{\otimes k}, C_*(\mte)),
$$ and   for simply connected $M$ there is a resulting chain equivalence (corollary 11 of \cite{cohenjones})
$$
f_* : C_*(\ltm) \to CH^*(C^*(\nee), C_*(\mte)),
$$ where  $CH^*(A,B)$ is the Hochschild cochain complex of a (differential graded) algebra $A$ with coefficients in the bimodule $B$.  
 Furthermore, as argued in the proof of theorem 13 of \cite{cohenjones}, $f_*: C_*(\ltm) \to CH^*(C^*(\nee), C_*(\mte))\cong CH^*(C_*(F(e)), C_*(\mte))$ takes the product coming from the ring spectrum structure of $\ltm$ to the cup product in the Hochschild cochain complex. 
 
 We therefore have a ring isomorphism $$\bh_*(LM) \cong H_*(\ltm) \xrightarrow{f_*} H^*(C^*(\nee), C_*(\mte))\cong H^*(C^*(M), C_*(\mtm)).$$ 
 
 Now on p. 794 of \cite{cohenjones} it was stated without proof that by ``S-duality" one can replace the coefficients $C_*(\mtm)$ by $C^*(M)$.  Proposition \ref{hoch} above supplies the justification for this assertion. Therefore we have a ring isomorphism, $$\bh_*(LM) \cong H^*(C^*(M), C^*(M))$$ as claimed.

\end{document}